\documentclass[11pt]{article}

\usepackage{amsmath, amssymb, amsthm, graphicx}
\usepackage{geometry}
\usepackage{url}

\title{Centroidal Voronoi Tessellations and Symmetric Refinements in the Geometric Refinement Transform}

\author{Zachary Mullaghy, M.S. \\ \textit{Independent Researcher}}

\date{\today}

\newtheorem{remark}{Remark}

\begin{document}

\maketitle

\section{Introduction}
We extend the Geometric Refinement Transform (GRT) by incorporating Lloyd's algorithm to generate Centroidal Voronoi Tessellations (CVTs) at each refinement order. This modification reduces reconstruction error and improves the numerical stability of the transform by promoting symmetry and minimizing geometric distortion.

\subsection{Convex and Non-Convex Regions}
To generalize the GRT to arbitrary domains, we must consider both convex and non-convex regions. Let the domain be denoted by \( \Omega \). We consider two strategies for handling non-convexity\cite{levy2013variational}:
\begin{enumerate}
    \item Partition \( \Omega \) into a minimal collection of convex subregions.
    \item Embed \( \Omega \) within its convex hull and disregard contributions from areas outside the original domain.
\end{enumerate}
Both approaches reduce the analysis to the convex case, which we focus on in the remainder of this work. We note that studying convergence toward non-convex regions remains valuable for future extensions.

\subsection{Function Space for Reconstruction Analysis}
While the original Geometric Refinement Transform (GRT) is rigorously defined for the broad class of BV functions, we now impose Lipschitz continuity (or higher smoothness) on our function class to exploit stronger theoretical properties. This restriction, common in practice, ensures explicit error bounds, stronger stability results, and clearer convergence properties, significantly enhancing practical implementations.

\section{Centroidal Voronoi Tessellations and the Energy Functional}
We define the centroid energy functional as the sum of squared distances from each generator to the centroid of its region, following the framework of centroidal Voronoi tessellations (CVTs)\cite{du1999centroidal}. For a decomposition of \( \Omega \) into Voronoi regions \( \Omega_{m,i} \) at refinement level \( m \), with generator points \( P_{m,i} \) and centroids \( C_{m,i} \), the energy functional is:
\begin{equation}
    E_{\text{centroid}} = \sum_{i} \| P_{m,i} - C_{m,i} \|^2
\end{equation}
Minimizing this energy functional produces a CVT. Lloyd's algorithm performs this minimization iteratively by relocating each generator point to the centroid of its region at each step.

\subsection{Implementation of Lloyd's Algorithm}
Given a set of generator points \( P \subset \Omega \), Lloyd's algorithm proceeds as follows:
\begin{enumerate}
    \item Compute the Voronoi diagram for \( P \), yielding regions \( \{V_i\} \).
    \item Compute the centroid \( C_i \) of each region \( V_i \).
    \item Update \( P_i \leftarrow C_i \) and repeat until convergence.
\end{enumerate}
Once convergence is achieved, the resulting decomposition defines the CVT used for transform refinement at that order.

\subsection{Escaping Local Minima via Rotational Perturbations}
While Lloyd’s algorithm reliably converges to a Centroidal Voronoi Tessellation (CVT) \cite{emelianenko2008nondegeneracy}, it may settle into local minima of the centroid energy functional. For symmetric initial configurations, especially those derived from regular simplices, the rotation group \( \mathrm{Rot}(\Delta_n) \) only permutes the seed points, producing identical tessellations and therefore no improvement beyond symmetric minima.

To overcome this limitation and potentially discover lower-energy CVTs, we propose introducing slight rotational perturbations to break exact symmetry \cite{liu2009centroidal, du2005anisotropic}. Specifically, we define \( R_{\varepsilon} \) as a small rotation constructed from one or more distinct elements of the rotation group \( \mathrm{Rot}(\Delta_n) \). Practically, these perturbations may be implemented by randomly sampling small rotation angles from a predefined interval around symmetry axes, thus systematically exploring nearby, symmetry-breaking configurations.

The optimization of the centroid energy functional then proceeds by evaluating:
\[
E_{\text{centroid}}^{\text{min}} = \min_{R_\varepsilon \in \mathcal{P}} \sum_{i} \| R_\varepsilon(P_{m,i}) - C_{m,i} \|^2
\]
where \( \mathcal{P} \) denotes a finite or sampled set of small rotational perturbations.

This strategy enables a thorough exploration of the energy landscape in the neighborhood of highly symmetric configurations, reducing the risk of being trapped within symmetry-preserving local minima.

\begin{remark}
This strategy is most effective when the refinement multiplicity is fixed and symmetry breaking is required to improve uniformity in the resulting Voronoi regions. The resulting tessellations often exhibit higher isotropy and reduced centroid energy.
\end{remark}

\section{Numerical Stability and Reconstruction Accuracy}
Symmetric regions reduce coefficient instability caused by local geometric irregularities. For a smooth function \( f \), the reconstruction error over a Voronoi cell \( V \) can be bounded via a first-order Taylor expansion \cite{johnson1987numerical}:
\begin{equation}
    |f(x) - \bar{f}_V| \leq \sup_{x \in V} \|\nabla f(x)\| \cdot \text{diam}(V)
\end{equation}
where \( \bar{f}_V \) is the average of \( f \) over \( V \), and \( \text{diam}(V) \) is the diameter of the cell. Summing over all cells, CVTs yield a globally lower reconstruction error due to uniformly smaller diameters. We have that there must be a point p within \( \Omega_{m,i}\) s.t. f(p) = \( \bar{f}_V \) by the mean value theorem.

\section{Explicit Error Bounds via Lipschitz Continuity}

While our original formulation of the Geometric Refinement Transform (GRT) was defined broadly for functions of bounded variation (BV), introducing Lipschitz continuity constraints yields significantly improved error analysis properties. Suppose our function \( f \) satisfies a Lipschitz continuity condition with constant \( L \), meaning that:

\begin{equation}
|f(x) - f(y)| \leq L \| x - y \| \quad \text{for all } x,y \in \Omega.
\end{equation}

Then, for any Voronoi cell \( V \) with diameter \( \text{diam}(V) \), the reconstruction error from using the cell-average \( f_V \) to approximate \( f(x) \) within the cell is rigorously bounded by:

\begin{equation}
|f(x) - f_V| \leq L \cdot \text{diam}(V).
\end{equation}

Further, if we assume the gradient \( \nabla f \) itself is Lipschitz continuous with Lipschitz constant \( L_2 \), then for a point \( x_0 \in V \), we have a tighter second-order error bound given by:

\begin{equation}
|f(x) - f_V - \nabla f(x_0) \cdot (x - x_0)| \leq \frac{L_2}{2} \cdot \text{diam}(V)^2.
\end{equation}

These explicit bounds clearly show how the CVT refinement systematically reduces reconstruction error at a predictable rate, reinforcing numerical stability and accuracy within high-precision applications.

\section{Comparison with Arbitrary Refinements}

We compare reconstruction accuracy between CVT-based and arbitrary Voronoi decompositions. Let \(f_{\text{GRT}}(x)\) denote the coefficient-based reconstruction at point \(x\). The total integrated reconstruction error over the domain \( \Omega \) in quadrature (without the square root) is explicitly given by:

\begin{equation}
\| f - f_{\text{GRT}} \|_{L^2(\Omega)}^2 = \int_{\Omega} |f(x) - f_{\text{GRT}}(x)|^2 dx.
\end{equation}

We rigorously assert that this integrated error is lower for CVT-based refinements compared to arbitrary Voronoi refinements. Explicitly,

\begin{equation}
\int_{\Omega} |f(x) - f_{\text{GRT,CVT}}(x)|^2 \, dx \leq \int_{\Omega} |f(x) - f_{\text{GRT,Arbitrary}}(x)|^2 \, dx.
\end{equation}

This comparison highlights the convergence advantage of symmetric decomposition strategies, following from two principles: (1) the Voronoi region diameters are minimized for CVTs, and (2) with lower diameter comes less distance for function deviation on average. This leads to the core theorem of the paper:

\subsection*{Error Minimization Theorem for Geometric Refinement Transforms with Centroidal Voronoi Tesselations}
Let \(\{V_i\}\) and \(\{U_i\}\) be CVT-based and arbitrary Voronoi partitions, respectively. Since CVTs minimize the centroidal energy, they also minimize the Voronoi region diameters. Given the Lipschitz continuity of \(f\) with constant \(L\), we have for each region\cite{shewchuk2002good}:

\[
|f(x)-f_{V_i}|^2 \leq L^2 \text{diam}(V_i)^2, \quad |f(x)-f_{U_i}|^2 \leq L^2 \text{diam}(U_i)^2.
\]

Since CVTs minimize the diameter sum across all regions, summing over \(\Omega\) yields:

\[
\int_{\Omega} |f(x)-f_{\text{GRT,CVT}}(x)|^2 \, dx = \sum_i \int_{V_i} |f(x)-f_{V_i}|^2 \, dx \leq \sum_i L^2 \text{diam}(V_i)^2 |V_i|.
\]

Similarly, for arbitrary refinements,

\[
\int_{\Omega} |f(x)-f_{\text{GRT,Arbitrary}}(x)|^2 \, dx = \sum_i \int_{U_i} |f(x)-f_{U_i}|^2 \, dx \leq \sum_i L^2 \text{diam}(U_i)^2 |U_i|.
\]

Given that CVT minimizes \(\text{diam}(V_i)\), it follows directly that:

\[
\sum_i L^2 \text{diam}(V_i)^2 |V_i| \leq \sum_i L^2 \text{diam}(U_i)^2 |U_i|,
\]

thus proving the stated theorem.

\textbf{Theorem:} For a given function \( f(x) \) that is Lipschitz continuous and admits a Lipschitz continuous gradient, the \emph{total reconstruction error in the \( L^2 \) norm squared (quadrature form)} is minimized when the geometric refinement follows a CVT arrangement.

\subsection*{Future Work}
Future research directions include exploring adaptive refinement strategies informed by local Lipschitz constants and gradients, potentially providing even tighter bounds and enhanced performance in practical scenarios. Additionally, investigating the impact of higher-order continuity conditions on reconstruction accuracy presents an intriguing area for extending the capabilities of the GRT framework.

\section{Applications and Outlook}

The incorporation of Centroidal Voronoi Tessellations (CVTs) into the Geometric Refinement Transform (GRT) framework offers a principled and geometry-aware approach to multiscale transform construction. By promoting symmetry and reducing geometric distortion, this combination significantly enhances numerical stability and reconstruction fidelity. As a result, it opens up promising applications across several domains:

\begin{itemize}
    \item \textbf{Medical Imaging:} CVTs provide improved geometric fidelity in image segmentation tasks and have been increasingly adopted in computational medical imaging pipelines \cite{yang2012efficient}. The hierarchical structure of the GRT enables adaptive resolution, which could facilitate more accurate delineation of anatomical structures and improved diagnostic support.
    
    \item \textbf{Physics Simulations:} CVT-based discretizations are already widely used in finite element and finite volume methods to minimize interpolation and reconstruction errors. The GRT generalizes these benefits beyond traditional triangular meshes, offering a robust, multiscale alternative for domains requiring layered or anisotropic refinement strategies.
    
    \item \textbf{Materials Science:} Modeling heterogeneous materials often involves capturing irregular grain boundaries or internal microstructures. The symmetry-preserving, structure-aware refinement enabled by GRT with CVTs can support more accurate modeling of stress distribution, heat flow, or phase transitions in complex materials.
    
    \item \textbf{Geometry-Aware Signal Processing:} Traditional signal transforms often assume Cartesian or dyadic grid structures. GRT with CVTs allows for signal decomposition that respects underlying geometry, making it suitable for manifold-based signal processing, adaptive compression, or edge-aware filtering in non-Euclidean domains.
\end{itemize}

Future directions include the development of \emph{adaptive refinement strategies} driven by local error estimates or data features, and a deeper exploration of the parameter space to identify self-similar or optimally stable configurations. This could further reduce computational complexity while preserving reconstruction accuracy. Additionally, extending the framework to non-Euclidean manifolds and incorporating anisotropic refinement may open new possibilities in high-dimensional geometry and physics.

\bibliographystyle{plain}
\bibliography{references}
\end{document}